\input amssym.def
\input amssym
\magnification=1200
\parindent0pt
\hsize=16 true cm \baselineskip=13  pt plus .2pt $ $

\def\Z{\Bbb Z}
\def\A{\Bbb A}
\def\S{\Bbb S}

\def\R{\Bbb R}

\centerline {\bf On finite  groups acting on acyclic low-dimensional manifolds}

\bigskip \bigskip

\centerline {Alessandra Guazzi*, Mattia Mecchia**  and  Bruno Zimmermann**}

\bigskip

\centerline {*SISSA} \centerline {Via Bonomea 256} \centerline {34136 Trieste,
Italy}

\bigskip

\centerline {**Universit\`a degli Studi di Trieste} \centerline {Dipartimento
di Matematica e Informatica} \centerline {34100 Trieste, Italy}

\bigskip \bigskip

Abstract.  We consider finite groups which admit a faithful, smooth action on
an acyclic manifold of dimension three, four or five (e.g. euclidean space).
Our first main result states that a finite group acting  on an acyclic 3- or
4-manifold is isomorphic to a subgroup of the orthogonal group ${\rm O}(3)$ or
${\rm O}(4)$, respectively. The analogue remains open in dimension five (where
it is not true for arbitrary continuous actions, however). We prove that  the
only finite nonabelian simple groups admitting a smooth action on an acyclic
5-manifold are the alternating groups $\A_5$ and  $\A_6$, and deduce from this
a short list of finite groups, closely related to the finite subgroups of
SO(5), which are the candidates for orientation-preserving actions on acyclic
5-manifolds.

\bigskip \bigskip

{\bf 1. Introduction}

\bigskip

All finite group actions considered in the present paper will be faithful and
smooth (or locally linear).

\medskip

By the recent geometrization of finite group actions on 3-manifolds, every
finite group action on the 3-sphere is conjugate to an orthogonal action; in
particular, the finite groups which occur are exactly the well-known finite
subgroups of the orthogonal group ${\rm O}(4)$. Finite groups acting on
arbitrary homology 3-spheres are considered in [MeZ1] and [Z]; here some other
finite groups occur and the situation is still not completely understood (see
[Z] and section 7).

\medskip

In dimension four, it is no longer true that a finite group action on the
4-sphere is conjugate to an orthogonal action (e.g. the Smith conjecture does
not remain true for the 4-sphere, that is the fixed point set of a periodic
diffeomorphism of $S^4$ may be a knotted 2-sphere). However it has been shown
in [MeZ2] and [MeZ3] that a finite group which admits an orientation-preserving
action on the 4-sphere, and more generally on any homology 4-sphere, is
isomorphic to a subgroup of the orthogonal group ${\rm SO}(5)$ (up to 2-fold
extensions in the case of solvable groups).

\medskip

In the present paper we consider finite groups acting on acyclic (compact or
non-compact) low-dimensional manifolds, i.e. manifolds with trivial reduced
integer homology (e.g. euclidean spaces).  Our first main result is the
following.

\bigskip

{\bf Theorem 1.}  {\sl  A finite group which admits a faithful, smooth action
on an acyclic 3- or 4-manifold is isomorphic to a subgroup of the orthogonal
group ${\rm O}(3)$ or ${\rm O}(4)$, respectively, and to a subgroup of  ${\rm
SO}(3)$ or ${\rm SO}(4)$ if the action is orientation-preserving.  In
particular, the only finite nonabelian simple group admitting such an action is
the alternating group $\A_5$.}

\bigskip

See [DV] for the finite subgroups of ${\rm O}(3)$ and ${\rm O}(4)$. Theorem 1
answers [E, Problem 11] on finite groups acting on euclidean space $\R^4$, for
the case of smooth actions.

\medskip

In the situation of Theorem 1, whereas each solvable group admitting an action
has a global fixed point, this does not remain true in general for nonsolvable
groups. As a classical example, the Poincar\'e homology 3-sphere admits an
action of $\A_5$ with a single fixed point, and the complement of this fixed
point is an acyclic 3-manifold with a fixed point free $\A_5$-action.  In
dimension five, we don't know the existence of a fixed point even for solvable
groups. On the other hand, as to nonsolvable groups, one of the main technical
problems in the proof of such classification results is to get hold of the
finite simple groups which may occur; in dimension five, the following is true.

\bigskip

{\bf Theorem 2.}  {\sl The only finite  nonabelian simple groups  admitting a
faithful, smooth action on an acyclic 5-manifold are the alternating groups
$\A_5$ and $\A_6$ (for actions of quasisimple groups there occurs, in addition,
the binary dodecahedral group $\A_5^*$).}

\bigskip

So these are exactly the simple (quasisimple) subgroups of the orthogonal group
${\rm SO}(5)$;  we recall that a quasisimple group is a perfect, central
extension of a simple group.  From Theorem 2 we deduce the following result for
arbitrary finite groups acting on acyclic 5-manifolds.

\bigskip

{\bf Theorem 3.}  {\sl Let $G$ be a finite group admitting a smooth, faithful,
orientation-preser-ving  action on  an acyclic 5-manifold. Then one of the
following cases occurs.

\smallskip

(i)  \hskip 2mm  $G$ is a subgroup of ${\rm SO}(5)$;

\smallskip

(ii) \hskip 1mm  $G$ contains a normal subgroup $N$ which is cyclic or a
central product of a cyclic group with $\A_5, \A_5^*$ or $\A_6$, and the factor
group $G/N$ is an elementary abelian 2-group  of rank at most four.}

\bigskip

See [MeZ3, Corollary 2] for a characterization of the finite subgroups of the
orthogonal group ${\rm SO}(5)$.  At present, we do not know an example of a
group which admits a faithful, smooth, orientation-preserving action on an
acyclic 5-manifold but is not isomorphic to a subgroup of ${\rm SO}(5)$;
however for continuous actions such examples in fact do exist. Specifically,
among the groups $G$ described in Theorem 3 (ii) there are the Milnor groups
$Q(8a,b,c)$ ([Mn]) which are extensions of a cyclic group by the elementary
abelian 2-group $(\Z_2)^2$.  Some of the Milnor groups admit a faithful,
continuous, orientation-preserving action on $\R^5$ but none of them is
isomorphic to a subgroup of ${\rm SO}(5)$; see Section 7.

\medskip

The paper is organized as follows. Section 2 contains some preliminary results
about finite group actions on acyclic manifolds. In Section 3 we present the
proof of Theorem 1 in the 3-dimensional case which is much shorter than the
4-dimensional (the main ingredient in dimension three is the Gorenstein-Walter
classification  of the finite simple groups with dihedral Sylow 2-subgroups, in
dimensions four and five this is replaced by the more involved
Gorenstein-Harada classification  of the finite simple groups of sectional
2-rank at most four); so this gives the reader a short approach to the basic
methods of the present paper, without the technical problems on finite simple
groups arising in dimensions four  and five.  In section 4 we prove  Theorem 2
concerning simple groups acting on acyclic 5-manifolds; this implies also the
analogue for simple groups acting on acyclic 4-manifolds, needed for the proof
of the 4-dimensional case of Theorem 1 given in section 5 (a direct proof in
dimension four would be somewhat shorter since the case of quasisimple groups
can be avoided; on the other hand, the main technical difficulties related to
the  Gorestein-Harada list remain the same, so we prefer to give the proof only
in dimension five). In section 6 we prove  Theorem 3, and in the last section
we discuss continuous actions of the  Milnor groups on acyclic 5-manifolds.

\bigskip  \bigskip

{\bf 2. Preliminary results.}

\bigskip

If $G$ is a finite group acting  on a $n$-manifold  the fixed point set of $G$
is  a submanifold and by Newman's Theorem   has dimension strictly smaller then
$n$ (see [B, Chapter III]).   If $G$ has  non-empty global fixed point set, the
group  $G$ leaves invariant  a tubular neighborhood of the fixed point set (see
[B, Chapter VI]). Hence, once we have found that $G$ fixes pointwise a
submanifold of dimension $d$,  we automatically get that it is a subgroup of
${\rm O}(n-d)$, and in particular of $\rm{SO}(n-d)$ if the action is
orientation-preserving.

Suppose that   $G$ is a $p$-group acting smoothly on  a $\Z_p$-acyclic
$n$-manifold (a $n$-manifold with trivial homology with coefficients in $\Z_p$,
the integers mod $p$). By Smith Theory (see [B, Chapter III, Section 5])  the
fixed point set of $G$  is again a $\Z_p$-acyclic manifold  (of even
codimension if $p$ is odd) and, in particular, it is nonempty; thus, $G$ is a
subgroup of ${\rm O}(n)$.
\medskip
By the above discussion we can state the following:

\bigskip

{\bf Lemma 1.} {\sl A finite $p$-group acting  on an acyclic $n$-manifold is
isomorphic to a subgroup of ${\rm O}(n)$, and in particular to a subgroup of
${\rm SO}(n)$ if the action is orientation-preserving}

\bigskip

We  note that every action of a finite group $G$ on an acyclic $1$- or
$2$-manifold has a global fixed point, so that, in particular, an analogous of
Theorem 1 also holds in the $1$- and $2$- dimensional case. The $1$-dimensional
case is obvious: an  acyclic 1-manifold (even $\Z_p$-acyclic)  is diffeomorphic
to $\Bbb{R}$, $[0,\infty)$ or to $[0,1]$, hence a finite group acting on it is
either trivial or generated by an orientation-reversing involution with one
fixed point. The 2-dimensional case is considered in the following lemma.

\bigskip

{\bf Lemma 2.} {\sl Any faithful and smooth action of  a finite group $G$ on an
acyclic $2$-manifold $X$ admits a global fixed point and $G$ is therefore
either dihedral with an orientation-reversing involution or cyclic.}

\bigskip

{\it Proof.} Suppose first that  $G$  is a nonabelian simple group; then the
action is orientation-preserving (since otherwise  $G$ would have a subgroup of
index two). Let $S$ be a Sylow 2-subgroup of $G$. By Lemma 1, the group $S$ is
a finite subgroup of ${\rm SO}(2)$, namely $S$ is cyclic; this is a
contradiction since a simple cannot have a cyclic Sylow 2-subgroup (see [Su2,
Corollary 2 of Theorem 2.2.10, p. 144]).

\medskip

Hence we can suppose that $G$ admits no nonabelian simple subgroups and, if $N$
is a minimal nontrivial normal subgroup in $G$, it is an elementary abelian
$p$-group by [Su1, Corollary 3 of  (2.4.14), p. 137]. Let $X^N$ be the
submanifold of points fixed by $N$; it has dimension at most $1$, is invariant
under the action of $G$ and $\Z_p$-acyclic. So either $X^N$ is  a point and we
are done, or $X^N$ is an acyclic $1$-manifold. In the latter case, let $T$ be
the normal subgroup of $G$ fixing all the points in $X^N$, the factor group
$G/T$ acts faithfully on $X^N$, with at least one global fixed point as noted
above, and also $G$  fixes that point. Hence $G$ is a finite subgroup of ${\rm
O}(2)$, and in particular of ${\rm SO}(2)$ if the action is
orientation-preserving.

\medskip

This concludes the proof of Lemma 1.

\bigskip

Note that every orientable $\Z_p$-acyclic 2-manifold  is an acyclic 2-manifold.
This can easily be seen, using simplicial homology and considering a cycle
$\alpha$ that is not an integral border (note that $\alpha$ can always be
chosen to have a connected and simple geometric realization), but that is a
border mod $p$ of a $2$-chain $\beta$ and reaching the contradiction that
$\alpha$ is an integral border. The orientability hypothesis cannot be omitted,
as it is needed to induce a coherent orientation on all $2$-simplices appearing
with nonzero coefficients in $\beta$. The projective plane is an example of
non-orientable $\Z_p$-acyclic 2-manifold  (with $p$ odd) that is not an acyclic
manifold.

\medskip

This remark, together with Lemma 2, leads us directly to a crucial lemma.

\bigskip

{\bf Lemma 3.} {\sl Let $G$ be a finite group acting on an acyclic $n$-manifold
$X$, with a nontrivial normal $p$-subgroup $N$. Suppose that one of the
following conditions holds:

\smallskip

(i)  \hskip 2mm the submanifold of points fixed by $N$ has dimension $d\leq 2$;

(ii)  \hskip 2mm  $n=3$;

(iii)  \hskip 1mm  $n=4$ and the action of $G$ is orientation-preserving;

(iv)  \hskip 1mm  $n=5$, the action of $G$ is orientation-preserving and $N$ is
not cyclic.

\smallskip

Then $G$ has at least one global fixed point. Hence $G$ acts orthogonally on
the boundary of some regular neighborhood of the fixed point and it is
isomorphic to a subgroup of ${\rm O}(n)$ (and in particular of ${\rm SO}(n)$ if
the action is orientation-preserving). }

\bigskip

{\it Proof.} Suppose that condition (i) holds.  If $d\leq 1$, we can proceed as
in the proof of Lemma 2.  Suppose therefore that $N$ fixes pointwise $X^N$, a
$\Z_p$-acyclic $2$-manifold. The 2-manifold $X^N$ is also orientable. In fact,
if $X^N$ is a $\Z_2$-acyclic manifold, it is orientable since otherwise the
first homology $H_1(X^N,\Z)$ would surject onto $\Z_2$; if $p$ is odd instead,
the order of $N$ is odd and $X^N$ is orientable by [B, Chapter IV, Theorem 2.1.
p.175]. Therefore, $X^N$ is an  acyclic $2$-manifold.

\medskip

Now, as in the proof of Lemma 2, we consider the normal subgroup $T$ of $G$
fixing each point of $X^N$. Then $G/T$ acts on $X^N$ and has a global fixed
point, hence also $G$ has a global fixed point. This concludes the proof in the
first case.

\medskip

Note that conditions (ii) and (iii) easily imply condition (i).

\medskip

If $n=5$ and the action of $G$ is orientation-preserving, the fixed point set
of $N$ may have also dimension three. In this case, by the discussion at the
beginning of the section, the subgroup $N$ is isomorphic to a finite subgroup
of  ${\rm SO(2)}$, namely it is  cyclic. This concludes the proof of the Lemma
3.

\bigskip  \bigskip

{\bf 3. Proof of Theorem 1 for acyclic 3-manifolds}

\bigskip

Let $G$ be a finite group with a smooth and  faithful action on an acyclic
3-manifold.

\bigskip

{\bf Proposition 1.} {\sl Suppose that $G$ is a nonabelian simple group acting
on an acyclic 3-manifold; then $G$ is isomorphic to the alternating group
$\A_5$.}

\bigskip

{\it Proof.} The action of $G$ is orientation-preserving as otherwise it should
have a subgroup of index two. By Lemma 1 a Sylow 2-subgroup of $G$  is
isomorphic to a subgroup of ${\rm SO}(3)$ and hence dihedral (since a Sylow
2-subgroup of a nonabelian simple group cannot be cyclic by [Su2, Corollary 2
of Theorem 2.2.10, p.144]). By the Gorenstein-Walter characterization of the
finite simple groups with dihedral Sylow 2-subgroups ([G1, Theorem 1.4.7],
[Su2, Theorem 6.8.6, p.505]),  $G$ is isomorphic to a linear fractional group
${\rm PSL}(2,q)$, for an odd prime power $q = p^n$, or to the alternating group
$\A_7$.

\medskip

Suppose first that $G$ is a linear fractional group ${\rm PSL}(2,q)$,  for an
odd prime power $q = p^n$. If $n > 1$ then ${\rm PSL}(2,q)$ has a noncyclic
elementary abelian $p$-subgroup $(\Z_p)^n$ (the subgroup represented by the
upper triangular matrices with both diagonal entries equal to one, isomorphic
to the additive group of the finite field with $q$ elements). By Lemma 1 the
group $(\Z_p)^n$ does not admit an action on an acyclic 3-manifold. Hence $n =
1$; now ${\rm PSL}(2,p)$ has a semidirect product  $\Z_p \rtimes \Z_{(p-1)/2}$
as a subgroup, with an effective action of $\Z_{(p-1)/2}$ (represented by the
diagonal matrices) on the normal subgroup $\Z_p$.  By Lemma 3, this is possible
only for $p = 5$, so we are left with the group ${\rm PSL}(2,5)$ isomorphic to
the alternating group $\A_5$.

\medskip

Finally, the alternating group  $\A_7$ has a subgroup $\Z_3 \times \Z_3$ which
is again excluded by Lemma 1.

\medskip

This completes the proof of Proposition 1.

\bigskip \bigskip

Recall that a {\it quasisimple} group is a perfect central extension of a
simple group, i.e. it is perfect and the factor group by its center is a
nonabelian simple group. A {\it semisimple} group is a central product of
quasisimple groups, i.e. the factor group by its center is a direct product of
nonabelian simple groups (see [Su2, chapter 6.6]). Any finite group $G$
contains a  unique maximal semisimple normal group $E(G)$ (which may be
trivial); the subgroup $E(G)$ is characteristic in $G$, and the quasisimple
factors of $E(G)$ are called the components of $G$.  The generalized Fitting
subgroup $F^*(G)$ of a group $G$ is defined as the (central) product of $E(G)$
and the Fitting subgroup $F(G)$ (the maximal nilpotent normal subgroup of $G$)
which is characteristic in  $G$. The generalized Fitting subgroup of a
nontrivial group is never trivial and its centralizer in $G$ coincides with its
center, i.e. $C_G(F^*(G))=Z(F^*(G))=Z(F(G))$ (see [Su2, Chapter 6.6, p.452]).

\bigskip  \bigskip

{\it Proof of Theorem 1 for 3-dimensional manifolds.}

\medskip

We divide the proof into two subcases.

\medskip

If the Fitting subgroup $F(G)$ is not trivial then $G$ has a non-trivial normal
$p$-subgroup, and by Lemma 3 the action of $G$ has a global fixed point.

\medskip

If the Fitting subgroup $F(G)$ is trivial, $E(G)$ coincides with generalized
Fitting subgroup $F^*(G)$, hence $C_G(E(G))=Z(E(G))=Z(F(G))$ is trivial.
Therefore, $G$ acts faithfully by conjugation on its normal subgroup $E(G)$,
namely, $G$ is a subgroup of $ {\rm Aut}(E(G))$, up to isomorphism. Also, by
definition of semisimple group, $E(G)\cong E(G)/Z(E(G))$ is a direct product of
simple groups acting on an acyclic $3$-manifold. Therefore, by Proposition 1,
$E(G)\cong(\A_5)^k$.

Suppose that $k\geq2$. Then $E(G)$ would contain an elementary abelian
$5$-subgroup of rank $2$, which is not a subgroup of ${\rm O}(3)$; a
contradiction. Therefore $E(G)$ is isomorphic to $\A_5$ and $G$ is a subgroup
of ${\rm Aut}(\A_5)\cong\S_5$, namely, either $G\cong \Bbb{S}_5$ or
$G\cong\A_5$. Suppose that $G\cong \S_5$. Then the subgroup generated by
$(4532),(12345)$ is a semidirect product $\Z_5\rtimes \Z_4$,  which is not a
subgroup of ${\rm O}(3)$. Thus, we have a contradiction by Lemma 3.

\medskip

We conclude that every finite group with trivial Fitting subgroup acting
smoothly on an acyclic $3$-manifold is isomorphic to $\A_5$,  a subgroup of
${\rm SO}(3)$, and this completes the proof of Theorem 1 for acyclic
3-manifolds.

\bigskip  \bigskip

{\bf 4. Proof of Theorem 2.}

\bigskip

To prove Theorem 2 we need some preliminary results.

\bigskip

{\bf Lemma 4.} {\sl  A finite group $G$ which admits a faithful,
orientation-preserving action on an acyclic 5-manifold has sectional 2-rank at
most four.}

\bigskip

{\it Proof.}  By Lemma 1, every $2$-subgroup of $G$ is in particular a subgroup
of ${\rm SO}(5)$ and it is therefore generated by at most 4 elements, (e.g. by
[MeZ2, Proposition 3.1.]).

\bigskip

{\bf Lemma 5.} {\sl  For a prime $p$ and a positive integer $r$, let $\Z_p
\rtimes \Z_r$ be a metacyclic group (semidirect product), with normal subgroup
$\Z_p$ and factor group $\Z_r$, which admits a faithful  and
orientation-preserving action on an acyclic 5-manifold. Then, by conjugation,
the square of each element of $\Z_r$ acts trivially or dihedrally on $\Z_p$.}

\bigskip

{\it Proof.}  If the fixed point set of $\Z_p$ is a 3-manifold, then $\Z_p$
locally acts as a group of rotations around it and an element in $\Z_r$
conjugates a rotation of minimal angle to a rotation of minimal angle, namely
it acts trivially or dihedrally on $\Z_p$. If instead $\Z_p$ fixes pointwise a
submanifold of dimension at most $2$, then, by Lemma 2, the group $\Z_p \rtimes
\Z_r$ is a subgroup of ${\rm SO}(5)$ and the claim follows from [MeZ2, Lemma
2.2].

\bigskip

We state also the following algebraic Lemma 6 that will be frequently used in
the proof of Theorem 2 and that is a simple consequence of [Su1, Theorem
2.9.18, p.257].

\bigskip

{\bf Lemma 6.} {\sl Let  $S$ be  a simple group; if $H$ is a simple subgroup of
$S$  then any central perfect extension  of $S$ contains a central perfect
extension of $H$.}

\bigskip  \bigskip

{\it Proof of Theorem 2.}

\bigskip

Let $G$ be a finite nonabelian quasisimple group acting on an acyclic
5-manifold. Since $G$ is perfect and it has no subgroup of index two, the
action of $G$ is orientation-preserving. By Lemma 4, the Sylow 2-subgroup $S$
has sectional 2-rank at most four. By the Gorenstein-Harada classification of
the simple groups of sectional 2-rank at most four ([G1, p.6], [Su2, Theorem
6.8.12, p.513]), the factor of $G$ by its center  is one of the groups in the
following list ($q$ denotes an odd prime power):

$${\rm PSL}(m,q), \;\; {\rm PSU}(m,q), \; \;  \; m \le 5,$$
$${\rm G}_2(q),\;\; ^3{\rm D}_4(q),\;\;  {\rm PSp}(4,q), \; \; ^2{\rm
G}_2(3^{2m+1}) \;\; (m \ge 1),$$
$${\rm PSL}(2,8),\;\; {\rm PSL(2,16)},\;\; {\rm PSL}(3,4),\;\; {\rm PSU}(3,4),\;\;
{\rm Sz}(8),$$
$$\A_m \;\; (7 \le m \le 11), \;\; {\rm M}_i \;\;(i \le 23),\;\; {\rm J}_i
\;\; (i\le 3),  \;\; {\rm McL}, \;\; {\rm Ly}.$$

\bigskip

In the following, we will exclude all the central perfect extensions of these
groups  except  $\A_5$, $\A_5^*$, and $\A_6$.

\medskip

We suppose first that $G$ is isomorphic to ${\rm SL}(2,p)$ or to ${\rm
PSL}(2,p)$, for   a prime $p \ge 5$. The group ${\rm SL}(2,p)$ has a metacyclic
subgroup $\Z_p \rtimes \Z_{p-1}$ (represented by all upper triangular
matrices): the normal subgroup $\Z_p$ consists of the matrices having both
entries on the diagonal equal to one and the subgroup  $\Z_{p-1}$ consists of
the diagonal matrices. The projection of $\Z_p \rtimes \Z_{p-1}$ to ${\rm
PSL}(2,p)$ is a metacyclic subgroup $\Z_p \rtimes \Z_{(p-1)/2}$,  and the
action  of $\Z_{(p-1)/2}$ on the normal subgroup $\Z_p$ is effective. By Lemma
5  we obtain that  $p=5$ and $G$ can be isomorphic to ${\rm PSL}(2,5)\cong
\A_5$ or to  ${\rm SL}(2,5)\cong \A_5^ *$ .

\medskip

Next we consider $G$ isomorphic to ${\rm PSL}(2,q)$ or to   ${\rm SL}(2,q)$ for
$q =p^n$ with $p$ odd prime number  and $n>1$.  In  ${\rm SL}(2,q)$ the
subgroup of upper triangular matrices is a semidirect product $(\Z_p)^n \rtimes
\Z_{q-1}$. The projection of the subgroup to ${\rm PSL}(2,q)$ is a semidirect
product  $(\Z_p)^n \rtimes \Z_{(q-1)/2}$, and the action of $\Z_{(q-1)/2}$ on
the normal subgroup  is effective.   In any case we can suppose that the group
contains a subgroup isomorphic to $(\Z_p)^n \rtimes \Z_r$ where   $r$ changes
according on the group we consider. The fixed point set of the subgroup
$(\Z_p)^n$ is an acyclic 1-manifold $M$; in fact, it cannot be  a
$\Z_p$-acyclic 3-manifold
 since  $(\Z_p)^n$ would  act faithfully on an orthogonal 1-sphere around
$M$. Now $(\Z_p)^n \rtimes \Z_r$ has a global fixed point, acts faithfully on
an orthogonal 3-sphere around $M$ in a fixed point and is a subgroup of the
orthogonal group O(4). This is possible only for  $n=2$ and $p=3$  (e.g.
$r=2^{k}$ with $k\leq 3$ by [MeZ1, Proposition 3 and 4]). The group ${\rm
SL}(2,9)$  can be excluded since  it contains a subgroup $ (\Z_3)^2\rtimes
\Z_{4}$ acting faithfully and orientation-preservingly on an orthogonal
3-sphere around $M$ but the  elements of order four in $\Z_{4}$ act dihedrally
on the subgroup $(\Z_3)^2$ and this is impossible (e.g.  by [MeZ1, Proposition
3]).  In this case $G$ is isomorphic to ${\rm PSL}(2,9)
 \cong  \A_6$.

\medskip

In general the unique central perfect extension of ${\rm PSL}(2,q)$ is ${\rm
SL}(2,q)$; the only exception is  ${\rm PSL}(2,9)$  that  has two other central
extensions one with center of order three and one with  center of order six
(see [Co, Table 5]). By [Co] we deduce that both these extensions do not
contain any element of order nine. In both cases the Sylow 3-subgroup of order
27 contains a normal subgroup isomorphic to $(\Z_3)^2$ and by the same argument
used  above, we can conclude that these groups cannot act on an acyclic
5-manifold.

\medskip

Now we consider  ${\rm PSL}(m,q)$, ${\rm SL}(m,q)$, ${\rm PSU}(m,q)$, ${\rm
SU}(m,q)$ with $q$ odd and $3\leq m \leq 5$ (for $m=2$ we have ${\rm
PSL}(2,q)\cong {\rm PSU}(2,q)$).  The groups  ${\rm PSL}(m,q)$ and ${\rm
SL}(m,q)$ contain a subgroup isomorphic  to ${\rm SL}(m-1,q)$ and the groups
${\rm PSU}(m,q)$ and ${\rm SU}(m,q)$ contain a subgroup isomorphic to ${\rm
SU}(m-1,q)$.  Using these subgroups we can eliminate inductively most part of
the groups. This process leaves   to consider case by case the simple groups
with $m=3$ and $q=3,5$ (these groups do not admit any central perfect extension
with non trivial center).   To eliminate  these groups  we use  [Co] to find
subgroups (metacyclic or simple) that we have already excluded.

\medskip

Also in this case in most part of cases the unique central extension of ${\rm
PSL}(m,q)$  (resp. ${\rm PSU}(m,q)$) is ${\rm SL}(m,q)$ (resp. ${\rm
SU}(m,q)$). We can have  an intermediate extension between ${\rm PSL}(4,q)$ and
${\rm SL}(4,q)$ and between  ${\rm PSU}(4,q)$ and ${\rm SU}(4,q)$ but the above
argument works again.  The only simple group of this type that admits further
central extensions is  ${\rm PSU}(4,3)$ (see [Co, Table 5]); in this case Lemma
6 and the inclusion  ${\rm PSL}(2,7) \subset {\A_7} \subset  {\rm PSU}(4,3)$
exclude directly  all the central perfect extensions.

\medskip

The groups ${\rm PSp}(4,q)$  contain a subgroup isomorphic to ${\rm PSL}(2,q)$.
This inclusion excludes automatically most part of the simple groups and their
central perfect extension (by Lemma 6); only  few  groups has to be checked
case by case. The group ${\rm PSp}(4,3)$ and its central extension contain  a
subgroup isomorphic to $(\Z_3)^3$ that can be excluded by the same argument
used for  ${\rm PSL}(2,p^n)$ with $n\geq 3$. The group  ${\rm PSp}(4,5)$
contains ${\rm PSL(2,25)}$ (see [Co]) while ${\rm PSp}(4,9)$ is excluded as it
has a subgroup isomorphic to  ${\rm PSp}(4,3)$.

\medskip

Up to central extension we have  $^3{\rm D}_4(q) \supset {\rm G}_2(q) \supset
{\rm PSL}(3,q)$ (see [St, Table 0A8], [GL, Table 4-1]). These inclusions and
Lemma 6 exclude $^3{\rm D}_4(q)$, $ {\rm G}_2(q)$ and their central extensions.

\medskip

The Ree groups $^2{\rm G}_2(3^{2m+1})$ have one conjugacy class of involutions,
the centralizer of an involution is $\Z_2^2 \times {\rm PSL}(2,3^{2m+1})$ ([G2,
p.164]) so for $m \ge 1$ they do not act (the group $^2{\rm G}_2(3)$ is not
simple).

\medskip

We consider now the  simple groups of Lie type and even characteristic. The
group  ${\rm PSL}(2,2^n)$ with $n=3,4$  contains a semidirect product $(\Z_2)^n
\ltimes \Z_r$. If ${\rm PSL}(2,2^n)$ acted on a acyclic 5-manifold, by Lemma 3
the semidirect product $(\Z_2)^n \ltimes \Z_r$ would be a subgroup of ${\rm
SO}(5)$ and this is impossible as the group $\Z_r$ acts transitively by
conjugation on $(\Z_2)^n$ (see  [MeZ3, Lemma 1]).

\medskip

The group ${\rm PSU}(3,4)$  contains  a subgroup isomorphic to $\Z_{13} \rtimes
\Z_3$ (see [Co]) and this group cannot act by Lemma 5. The group ${\rm
PSU}(3,4)$ does not admit any  central extension with non trivial center.

\medskip

The group ${\rm PSL}(3,4)$  contains a subgroup isomorphic to ${\rm PSL}(2,7)$
that we have already excluded. To eliminate the central extensions we apply
Lemma 6.

\medskip

Concerning the Suzuki group ${\rm Sz}(8)$, its Sylow 2-subgroup has order 64
with a normal subgroup $(\Z_2)^3$ and it has a unique conjugacy class of
involutions (see [Co]). If ${\rm Sz}(8)$ admitted an action, the subgroup
$(\Z_2)^3$ should have a global fixed point and  should act on the $4$-sphere
that is the boundary of some regular neighborhood of the fixed point. This is
impossible as the involutions in $(\Z_2)^3$ are all conjugated (see [MeZ3,
Lemma 1]). Suppose now that  $G$ is a  central prefect extension of ${\rm
Sz}(8)$; the center of  $G$  is an elementary abelian 2-group of rank one or
two. In any case the center fixes pointwise a $\Z_2$-acyclic manifold  $M$  of
dimension at most three. Since ${\rm Sz}(8)$ is simple, the normal subgroup of
the elements of $G$ leaving invariant each point of $M$  coincides with the
center of $G$. Hence, the quotient of $G$ by its center, and thus  $(\Z_2)^3$,
act  faithfully and orientation-preservengly on $M$,  which is impossible by
Lemma 1.

\medskip

The alternating group $\A_7$ has a subgroup ${\rm PSL}(2,7)$ which excludes all
alternating groups $\A_n$ for $n \ge 7$. For any of the remaining  simple
groups is possible to find a simple subgroup already excluded. We will not give
further details and refer to [Co] and its references for the maximal subgroups.

\medskip

This concludes the proof of Theorem 2.

\bigskip  \bigskip

{\bf 5. Proof of Theorem 1 for acyclic 4-manifolds}

\bigskip

In dimension $4$, Lemma 3 still applies for groups with nontrivial Fitting
subgroup if we suppose that the action is orientation-preserving. Hence, let us
first suppose that $G$ is a finite group admitting a smooth, faithful and
orientation-preserving action on an acyclic 4-manifold. We will then extend our
results to any smooth and faithful action.

\bigskip

{\bf Proposition 2.} {\sl A finite group $G$ which admits a smooth and
orientation-preserving action on an acyclic 4-manifold $X$ is isomorphic to a
subgroup of ${\rm SO}(4)$. In particular if  $F(G)$ is trivial, $G$ is
isomorphic to the alternating group $\A_5$.}

\bigskip

{\it Proof.} Suppose first that $G$ is simple and nonabelian. Then $G$ acts
also on $X\times \R$, an acyclic $5$-manifold and, by Theorem 2, the group $G$
is isomorphic to either $\A_5$ or $\A_6$. But $\A_6\cong {\rm PSL}(2,9)$
contains a solvable subgroup of the form $(\Z_3)^2\rtimes \Z_4$, which is not a
subgroup of ${\rm SO}(4)$ (it is shown in [MeZ1, proof of Theorem 2] that it
does not even act on a homology $3$-sphere); a contradiction, by Lemma 3.
Hence, $\A_5$ is the only nonabelian simple group which can act on an acyclic
$4$-manifold.

\medskip

Suppose now that the Fitting subgroup of $G$ is trivial; then $E(G)$ coincides
with the generalized Fitting $F^*(G)$, and hence $C_G(E(G))=Z(E(G))=Z(F^*(G))$.
Therefore $G$ is isomorphic to a subgroup of ${\rm Aut}(E(G))$ and
$E(G)=E(G)/Z(E(G))$ is a (nontrivial) direct sum of simple nonabelian groups,
acting on an acyclic $4$-manifold; thus $E(G)\cong(\A_5)^k$. For $k\geq 2$,
$(\A_5)^k$ contains an elementary abelian $2$-subgroup of rank at least $4$,
which is not a subgroup of ${\rm SO}(4)$ (e.g. by [MeZ1, Proposition 4]).
Hence, $\A_5\cong E(G)\subset G\subset {\rm Aut(E(G))}\cong \S_5$ and we deduce
that either $G\cong \A_5$ or $G\cong \S_5$. Finally, $\S_5$ contains a subgroup
isomorphic to $\Z_5\rtimes \Z_4$ (with faithful action of $\Z_4$ on $\Z_5$)
which is not a subgroup of ${\rm SO}(4)$  (by [Z, Proposition 3] it does not
even act on a homology $3$-sphere). We conclude that $G\cong\A_5$.

\medskip

If $F(G)$ is not trivial we can apply Lemma 3 and the thesis follows.

\medskip

This concludes the proof of Proposition 2.

\bigskip  \bigskip

{\it Proof of Theorem 1 for acyclic $4$-manifolds.}

\bigskip

Let  $G$ be a group acting on a acyclic 4-manifold $X$. We divide the proof
into two subcases.

\medskip

Suppose that the Fitting subgroup $F(G)$ is not trivial. If $G$ contains a
nontrivial normal $p$-subgroup $P$ such that the submanifold  of points fixed
by $P$ has dimension at most $2$, then Lemma 3 applies and the claim is proven.
Otherwise, $F(G)$ being nontrivial, $G$ admits a normal $p$-subgroup $P$ which
fixes pointwise a $3$-submanifold and is therefore generated by an
orientation-reversing involution $t$.

\medskip

Let $G_0$   be the index two subgroup of orientation-preserving elements in
$G$. As $t$ is an orientation-reversing element, $G_0\cap P=1$ and both $P$ and
$G_0$ are normal in $G$. We obtain that   $G=G_0\times P\cong G_0\times \Z_2$.
Note $F(G_0)$ is trivial, otherwise, $G_0$ (and hence $G\cong G_0\times P$)
would admit a normal $p$-subgroup which, acting orientation-preserving on $X$,
fixes a submanifold of dimension at most $2$; a contradiction. Therefore $G_0$
is a group acting orientation-preservingly on an acyclic $4$-manifold with
trivial Fitting subgroup.

\medskip

If $G_0$ is trivial, $G\cong \Z_2$, otherwise $G_0\cong\A_5$, by Proposition 2.
In the latter case, $G\cong \A_5\times\Z_2$, which is a subgroup of ${\rm
O}(3)$ and hence of ${\rm O}(4)$.

\medskip

If instead $F(G)$ is trivial, $E(G)=F^*(G)$ is not trivial. Being $E(G)$
semisimple, in particular $E(G)$ is perfect and hence its action on $X$ is
orientation-preserving; also, $F(E(G))$ is trivial, thus $E(G)\cong\A_5$, by
Proposition 2. Since the centralizer $C_G(E(G))=C_G(F^*(G))=Z(F^*(G))= Z(F(G))$
is trivial, $G$ is isomorphic to a subgroup of ${\rm Aut}(E(G))\cong \S_5$
containing $\A_5$. Therefore either $G\cong \A_5$ or $G=\S_5$, both subgroups
of ${\rm O}(4)$.

\medskip

For an orthogonal action of $\S_5$ on $\R^4$, just consider the action of
$\S_5$ on $\R^5$ permuting the standard orthonormal base  and restrict the
action to the hyperplane described by the equation $x_1+\ldots + x_5=0$, which
is invariant under the action of $\S_5$ on $\R^5$.

\medskip

This concludes the proof of Theorem 1.

\bigskip  \bigskip
\vfill \eject

{\bf 6. Proof of  Theorem 3}

\bigskip

Suppose that a group $G$ acts preserving the orientation on an acyclic
$5$-manifold and that it has nontrivial Fitting subgroup $F(G)$. By Lemma 3
either $G$ has a global fixed point and is a subgroup of ${\rm SO}(5)$,  or the
Fitting subgroup  $F(G)$ is a direct product of cyclic groups of coprime orders
fixing  $3$-manifolds. In the latter case, in particular, $F(G)$ is cyclic.

\medskip

This implies that if an abelian (or even nilpotent) group acts
orientation-preservingly and with no global fixed point on an acyclic
$5$-manifold, it is cyclic. As for cyclic groups which are not $p$-groups, it
is still an open question whether or not they can act with no global fixed
points, depending on how many primes divide their order (see [HKMS]). Note that
if the acyclic manifold is homeomorphic to a closed disk this is not possible,
by Brower's Theorem.

\bigskip

{\bf Proposition 3.} {\sl Suppose that $G$ is a (nontrivial) semisimple group
acting on an acyclic 5-manifold; then $G$ is isomorphic to one of the following
groups:
$$\A_5,\;\;  \A_6, \;\;   \A_5^*, \;\;   \A_5^*\times_{\Z_2} \A_5^*.$$
In particular, a semisimple group can act on an acyclic $5$-manifold if and
only if it is a subgroup of ${\rm SO}(5)$. }

\bigskip

{\it Proof.} Since $G$ is perfect and nonabelian the action is
orientation-preserving.  By Theorem 2, the quasisimple components of $G$ are
isomorphic to   $\A_5,\,\A_6$ or $\A_5^*$. Since by Lemma 4 the sectional
2-rank of $G$ is at most four, $G$ is the central product of at most two of
these quasisimple groups, so it remains to analyze the case of groups with two
quasisimple components. Suppose first that $G=Q_1\times Q_2$, where $Q_1$ is
isomorphic either to $\A_5$ or to $\A_6$. In this case $G$ has a subgroup
isomorphic to $\A_4 \times Q_2$. The subgroup $\A_4 \times Q_2$ contains a
normal elementary 2-group of rank two; by Lemma 3 the group $\A_4 \times Q_2$
should be  a subgroup of ${\rm SO}(5)$ which is not the case (e.g. by [MeZ2]).

\medskip

Next, the group $\A_5^*\times \A_5^*$ has a normal elementary 2-subgroup of
rank two and again is not a subgroup of SO(5).  So the only remaining
semisimple group with two components is $\A_5^*\times_{\Z_2} \A_5^*$.

\medskip

This finishes the proof of Proposition 3.

\bigskip  \bigskip

{\it Proof of Theorem 3.}

\bigskip

By Lemma 3 we can suppose that either the Fitting subgroup $F(G)$ is trivial,
or $F(G)$ is cyclic and the fixed point set of each $p$-subgroup of $F$ is a
3-dimensional manifold (since otherwise $G$ is a subgroup of ${\rm SO}(5)$ and
we are done). In the latter case   each $p$-subgroup of $F(G)$ acts as a
rotation group around its 3-dimensional fixed point set. Each element of $G$
acts by conjugation on each $p$-subgroup of $F(G)$; this action may be either
trivial or dihedral since a rotation of minimal angle around the fixed point
set is conjugated to a rotation of minimal angle. In any case the square of
each element of $G$  acts by conjugation trivially on $F(G)$.

\medskip

Suppose first that   $E(G)$ is trivial; $F(G)$  coincides with the generalized
Fitting subgroup $F^*(G)$ of $G$ (the product of the Fitting subgroup $F$(G)
and the maximal normal semisimple subgroup $E(G)$), and $F^*(G) = F(G)$
contains its centralizer in $G$ ([Su2, Theorem 6.6.11, p.452]). By the
preceding paragraph, this implies that each element in $G/F(G)$ has order at
most two, so $G/F(G)$ is an elementary abelian 2-group, of rank at most four by
Lemma 4.

\medskip

Suppose now that  $E(G)$ is nontrivial, and hence isomorphic to one of the
groups in Proposition 3.  By [Su2, Theorem 6.6.11, p.452], the factor group  of
$G$ by the center of $F(G)$ is isomorphic to a subgroup of the automorphism
group ${\rm Aut}(F^*(G))$ of the generalized Fitting subgroup. In our
situation, $F(G)$ is  cyclic and hence $G/F(G)$ is isomorphic to a subgroup of
${\rm Aut}(F^*(G))$ which in turn is isomorphic to a subgroup of ${\rm
Aut}(F(G)) \times {\rm Aut}(E(G))$ (since $F(G)$ and $E(G)$ are characteristic
in $G$).

\medskip

If $E(G)$ is isomorphic to $\A_5$ or $\A_5^*$, the outer automorphism of $E(G)$
has order two; if $E(G)$ is isomorphic to $\A_6$ the outer automorphism group
is elementary abelian of order four (see [Co]). We have supposed  that the
square of each element in $G$ acts trivially on $F(G)$ and, if $E(G)$ has only
one component, such a square acts as an inner automorphism on $E(G)$; this
implies easily that the square of each element in $G$ is contained in $F^*(G)$,
and hence the factor group $G/F^*(G)$ is an elementary abelian 2-group, of rank
at most four.

\medskip

Suppose that $E(G)$ is isomorphic to  $\A_5^*  \times_{\Z_2} \A_5^*$; then the
center $\Z_2$ of $E(G)$ is normal in $G$. By Lemma 3, we can assume that the
fixed point set of this normal subgroup $\Z_2$ has dimension three and hence,
by Smith fixed point theory, is a $\Z_2$-acyclic 3-manifold $M$. The factor
group $E(G)/\Z_2$ is isomorphic to $\A_5 \times \A_5$ and admits a faithful,
orientation-preserving action on $M$. Now $\A_5  \times \A_5 $ has a subgroup
$(\Z_2)^2 \times (\Z_2)^2 = (\Z_2)^4$; however, again by Smith fixed point
theory, the group $(\Z_2)^4$ does not admit a faithful, orientation-preserving
action on a $\Z_2$-homology 3-sphere.  So, if $E(G) \cong \A_5^*  \times_{\Z_2}
\A_5^*$, we have shown that $G$ is a subgroup of ${\rm SO}(5)$.

\medskip

This finishes the proof of Theorem 3.

\bigskip \bigskip

{\bf 7. The Milnor groups $Q(8a,b,c)$}

\bigskip

It is observed in [Mn] that the groups $Q(8a,b,c)$ have periodic cohomology of
period four but do not admit faithful, linear, free actions on $S^3$ (in fact,
they are not isomorphic to subgroups of O(4)).  We will assume in the following
that $a > b > c \ge 1$ are odd coprime integers; then $Q(8a,b,c)$  is a
semidirect product $\Z_{abc} \ltimes Q_8$ of a normal cyclic subgroup  $\Z_a
\times \Z_b \times \Z_c  \cong \Z_{abc}$ by the quaternion group $Q_8 = \{\pm
1,\pm i,\pm j,\pm k\}$ of order eight, where $i,j$ and $k$ act trivially on
$\Z_a, \Z_b$ and $\Z_c$, respectively, and in a dihedral way on the other two.
Note that $Q(8a,b,c)$ has also a normal subgroup $\Z_{2abc}$, with factor group
the elementary abelian 2-group $\Z_2^2$, so it is one of the groups described
in Theorem 3.

\medskip

It has been shown by Milgram ([Mg]; see also the comments in [K, p.173, Update
A to Problem 3.37]) that some of the groups $Q(8a,b,c)$ admit a faithful, free
action on a homology 3-sphere; let $Q$ be one of these groups which admits such
an action on a homology 3-sphere $M$.  By the double suspension theorem (see
e.g. [Ca]), the double suspension $M * S^1$ of $M$ (the join of $M$ with the
1-sphere) is homeomorphic to $S^5$. Letting $Q$ act trivially on $S^1$, the
actions of $Q$ on $M$ and $S^1$ induce a faithful, continuous,
orientation-preserving action of $Q$ on $S^5$ with fixed point set $S^1$, and
hence also on $\R^5$ (the complement of a fixed point). Now it is not difficult
to show that none of the groups $Q(8a,b,c)$ is isomorphic to a subgroup of
SO(5). At present, we do not know if some Milnor group $Q(8a,b,c)$ admits a
faithful, smooth, orientation-preserving action on an acyclic 5-manifold.

\bigskip   \bigskip

\centerline {\bf References}

\bigskip

\item {[B]} G.Bredon, {\it Introduction to Compact Transformation Groups.}Academic
Press, New York 1972

\item {[Ca]}  J.W.Cannon, {\it  The recognition problem: what is a topological
manifold.} Bull. Amer. Math. Soc. 84, 832-866  (1978)

\item {[Co]} J.H.Conway, R.T.Curtis, S.P.Norton, R.A.Parker, R.A.Wilson, {\it Atlas
of Finite Groups.} Oxford University Press 1985

\item {[DV]}  P.Du Val, {\it  Homographies, Quaternions and Rotations.} Oxford
Math. Monographs, Oxford University Press 1964

\item {[E]} A.L. Edmonds, {\it A survey of group actions on 4-manifolds.}
arXiv:0907.0454

\item {[G1]} D.Gorenstein, {\it The Classification of Finite Simple Groups.}
Plenum Press, New York 1983

\item {[G2]} D.Gorenstein, {\it Finite Simple Groups: An Introduction to  Their
Classification.}    Plenum Press, New York 1982

\item {[GL]} D.Gorenstein, R.Lyons, {\it The local structure of finite group of
characteristic 2 type.} Memoirs Amer. Math. Soc. 276 (vol. 42), 1-731  (1983)

\item {[HKMS]} R.Haynes, S.Kwasik, J.Mast, R.Schultz,
{\it  Periodic maps on ${\bf R}^7$ without fixed points.} Math. Proc. Cambridge
Philos. Soc.  132, 131-136 (2002).

\item {[K]} R.Kirby, {\it Problems in low-dimensional topology.} Geometric
Topology, AMS/IP Studies in Advanced Mathematics Volume 2 Part 2 (W.H.Kazez,
editor), 35- 358  (1997)

\item {[MeZ1]} M.Mecchia, B.Zimmermann, {\it On finite groups acting on $\Bbb
Z_2$-homology 3-spheres.}  Math. Z. 248,   675-693  (2004)

\item {[MeZ2]} M.Mecchia, B.Zimmermann, {\it On finite simple and nonsolvable
groups acting on homology 4-spheres.}  Top. Appl. 153, 2933-2942 (2006)

\item {[MeZ3]} M.Mecchia, B.Zimmermann, {\it On finite groups acting on homology
4-spheres and finite subgroups of ${\rm SO}(5)$.}  arXiv:1001.3976

\item {[Mg]} R.J.Milgram, {\it Evaluating the Swan finiteness obstruction for
finite groups.} Algebraic and Geometric Topology. Lecture Notes in Math. 1126
(Springer 1985), 127-158

\item {[Mn]} J.Milnor, {\it Groups which act on $S^n$ without fixed points.} Amer.
J. Math. 79, 623-630  (1957)

\item {[St]} E.Stensholt, {\it Certain embeddings among finite groups of Lie
type.} J. Algebra 53, 136-187  (1978)

\item {[Su1]} M.Suzuki, {\it Group Theory II.}  Springer-Verlag 1982

\item {[Su2]} M.Suzuki, {\it Group Theory II.}  Springer-Verlag 1982

\item {[Z]} B.Zimmermann, {\it On the classification of finite groups acting on
homology 3-spheres.}  Pac. J. Math. 217, 387-395 (2004)

\bye